\documentclass[12pt]{article}
\usepackage{amssymb, amsmath, latexsym, a4}
\usepackage[latin1]{inputenc}
\usepackage{graphicx}
\usepackage[all]{xy}
\usepackage[usenames]{color}

\newcommand{\rdg}{\hfill $\Box $}

\newtheorem{De}{Definition}[section]
\newtheorem{Th}[De]{Theorem}
\newtheorem{Pro}[De]{Proposition}
\newtheorem{Le}[De]{Lemma}
\newtheorem{Co}[De]{Corollary}
\newtheorem{Rem}[De]{Remark}
\newtheorem{Ex}[De]{Example}

\newcommand{\Image}{{\sf Im}}
\newcommand{\Ker}{{\sf Ker}}
\newcommand{\HL}{{\sf HL}}
\newcommand{\Lie}{\ensuremath{\mathsf{Lie}}}

\newcommand{\Lieh}{\ensuremath{\mathfrak{h}}}
\newcommand{\Lieg}{\ensuremath{\mathfrak{g}}}
\newcommand{\Lieq}{\ensuremath{\mathfrak{q}}}
\newcommand{\Liea}{\ensuremath{\mathfrak{a}}}
\newcommand{\Lieb}{\ensuremath{\mathfrak{b}}}
\newcommand{\Liem}{\ensuremath{\mathfrak{m}}}
\newcommand{\Lien}{\ensuremath{\mathfrak{n}}}
\newcommand{\Lief}{\ensuremath{\mathfrak{f}}}
\newcommand{\Lies}{\ensuremath{\mathfrak{s}}}
\newcommand{\Lier}{\ensuremath{\mathfrak{r}}}

\newcommand{\Liei}{\ensuremath{\mathfrak{i}}}
\newcommand{\Leib}{\ensuremath{\mathsf{Leib}}}

\newbox\pullbackbox
\setbox\pullbackbox=\hbox{\xy 0;<1mm,0mm>: \POS(4,0)\ar@{-} (0,0) \ar@{-} (4,4)
\endxy}

\newdir{>>}{{}*!/3.5pt/:(1,-.2)@^{>}*!/3.5pt/:(1,+.2)@_{>}*!/7pt/:(1,-.2)@^{>}*!/7pt/:(1,+.2)@_{>}}
\newdir{ >>}{{}*!/8pt/@{|}*!/3.5pt/:(1,-.2)@^{>}*!/3.5pt/:(1,+.2)@_{>}}
\newdir{ |>}{{}*!/-3.5pt/@{|}*!/-8pt/:(1,-.2)@^{>}*!/-8pt/:(1,+.2)@_{>}}
\newdir{ >}{{}*!/-8pt/@{>}}
\newdir{>}{{}*:(1,-.2)@^{>}*:(1,+.2)@_{>}}
\newdir{<}{{}*:(1,+.2)@^{<}*:(1,-.2)@_{<}}

\newcommand{\q}{\frak q}

   \newcommand{\eq}{\frak q} 
  \newcommand{\eh}{\frak h}

\begin{document}

\centerline{\bf  ON $\Lie$-ISOCLINIC LEIBNIZ ALGEBRAS}

\bigskip
\bigskip
\centerline{\bf G. R. Biyogmam$^{(1)}$ and J. M. Casas$^{(2)}$}

\bigskip
\centerline{$^{(1)}$ Department of Mathematics, Southwestern Oklahoma State University}
\centerline{Weatherford, OK 73096, USA}
\centerline{E-mail address: {guy.biyogmam@swosu.edu}}

\medskip

\centerline{$^{(2)}$Dpto. Matemática Aplicada, Universidad de Vigo,  E. E. Forestal}
\centerline{Campus Universitario A Xunqueira, 36005 Pontevedra, Spain}
\centerline{ {E-mail address}: jmcasas@uvigo.es}

\date{}

\bigskip \bigskip

{\bf Abstract:}  In this paper we study the notion of isoclinism on \Lie-central extensions of Leibniz algebras, this yields to introduce the concept of \Lie-isoclinic Leibniz algebras. We provide several equivalent conditions under which Leibniz algebras are \Lie-isoclinic. We also define the concept of Schur \Lie-multiplier  and  analyze its connection with \Lie-isoclinism.

\bigskip

\bigskip



\section{Introduction}

The concept of isoclinism goes back to Philip Hall \cite{PH} in 1940 in his attempt to classify $p$-groups using an equivalence relation weaker than the notion of isomorphism. This concept was extended to Lie algebras by K. Moneyhun in \cite{KM}, and later by several other authors \cite{Mog, Moh, Sal1, Sal2, Mohs}.  Our main purpose on this paper is to extend the study of isoclinism on non Lie-Leibniz algebras, relying on the anti-symmetry of their Leibniz bracket. Nevertheless, the generalization of isoclinism concepts from Lie to Leibniz algebras seems to be a mere adaptation of techniques and results, so our approach focuses in the relative notions of isoclinism.

Relative means that we are approaching the concepts relative to the Liezation functor $(-)_{\Lie}: {\sf Leib} \to {\sf Lie}$ which assigns to a Leibniz algebra $\Lieq$ the Lie algebra $\Lieq_{\Lie} = \Lieq/<\{[x,x]: x \in \Lieq\}>$. With other words, from the categorical theory of central extensions relative to a chosen subcategory of a base category given in \cite{JK}, the concepts of central extension, perfect object and commutator in a semi-abelian category relative to a Birkhoff subcategory were approached in \cite{CVDL1}. In particular, in \cite{CK} was approached the case of central extensions of Leibniz algebras with respect to its Birkhoff subcategory of Lie algebras, that is, with respect to the Liezation functor. The absolute case in this theory corresponds to the abelianization functor which assigns to a Leibniz (Lie) algebra the abelian Leibniz (Lie) algebra with the same underlying vector space, but with trivial bracket operation. This absolute case corresponds to the classical case described at the beginning.

In this paper we continue the study of central extensions relative to the Liezation functor, so called $\Lie$-central extensions, initiated in \cite{CK}. In concrete, we study isoclinism of $\Lie$-central extensions of Leibniz algebras. To do so, we organize
the paper as follows, in Section \ref{preliminaries} we recall the $\Lie$-notions defined in \cite{CK}. In Section \ref{Lie iso Lb alg}, we define the notion of $\Lie$-isoclinism between two $\Lie$-central extensions of Leibniz algebras and provide necessary and sufficient conditions under which two Leibniz algebras are $\Lie$-isoclinic.  Also in \cite{CK}, the authors defined the $\Lie$-homology of Leibniz algebras which enables us to analyze the connection between $\Lie$-isoclinism and the Schur $\Lie$-multiplier in Section \ref{Schur}. In this section, we provide another characterization of $\Lie$-isoclinic Leibniz algebras by means of the $\Lie$-homology of the isomorphisms between their $\Lie$-centers  and their $\Lie$-commutators. This is used to determine the group of $\Lie$-autoclinisms of any given Leibniz algebra. Finally we prove that under certain conditions, all $\Lie$-stem covers of a given Leibniz algebra are mutually $\Lie$-isoclinic.


\section{Preliminary results on Leibniz algebras} \label{preliminaries}
We fix $\mathbb{K}$ as a ground field such that $\frac{1}{2} \in \mathbb{K}$. All vector spaces and tensor products are considered over $\mathbb{K}$.

A \emph{Leibniz algebra} \cite{Lo 1, Lo 2} is a vector space $\eq$  equipped with a bilinear map $[-,-] : \Lieq \otimes \Lieq \to \Lieq$, usually called the \emph{Leibniz bracket} of $\eq$,  satisfying the \emph{Leibniz identity}:
\[
 [x,[y,z]]= [[x,y],z]-[[x,z],y], \quad x, y, z \in \Lieq.
\]
Leibniz algebras form a semi-abelian category \cite{CVDL, JMT}, denoted by {\Leib}, whose morphisms are linear maps that preserve the Leibniz bracket.

 A subalgebra ${\eh}$ of a Leibniz algebra ${\Lieq}$ is said to be \emph{left (resp. right) ideal} of ${\Lieq}$ if $ [h,q]\in {\eh}$  (resp.  $ [q,h]\in {\eh}$), for all $h \in {\eh}$, $q \in {\Lieq}$. If ${\eh}$ is both
left and right ideal, then ${\eh}$ is called \emph{two-sided ideal} of ${\Lieq}$. In this case $\Lieq/\Lieh$ naturally inherits a Leibniz algebra structure.

For a Leibniz algebra ${\Lieq}$, we denote by ${\Lieq}^{\rm ann}$ the subspace of ${\Lieq}$ spanned by all elements of the form $[x,x]$, $x \in \Lieq$.

Given a Leibniz algebra $\Lieq$, it is clear that the quotient ${\Lieq}_ {_{\rm Lie}}=\Lieq/{\Lieq}^{\rm ann}$ is a Lie algebra. This defines the so-called  \emph{Liezation functor} $(-)_{\Lie} : {\Leib} \to {\Lie}$, which assigns to a Leibniz algebra $\Lieq$ the Lie algebra ${\Lieq}_{_{\rm Lie}}$. Moreover,
the canonical epimorphism  ${\Lieq} \twoheadrightarrow {\Lieq}_ {_{\rm Lie}}$ is universal among all homomorphisms from $\Lieq$ to a Lie algebra, implying that the Liezation functor is left adjoint to the inclusion functor $ {\Lie} \hookrightarrow {\Leib} $.

The following notions come from \cite{CK}.

Let {\Lieq} be a Leibniz algebra and  ${\Liem}$, ${\Lien}$ be two-sided ideals of  ${\Lieq}$. The \emph{$\Lie$-commutator} of  ${\Liem}$ and  ${\Lien}$ is the two-sided ideal  of $\Lieq$
\[
[\Liem,\Lien]_{\Lie}= \langle \{[m,n]+[n,m], m \in \Liem, n \in \Lien \}\rangle.
\]

The \emph{$\Lie$-center} of the Leibniz algebra $\Lieq$ is the two-sided-ideal
\[
Z_{\Lie}(\Lieq) =  \{ z\in \Lieq\,|\,\text{$[q,z]+[z,q]=0$ for all $q\in \Lieq$}\}.
\]

An extension of Leibniz algebras $0 \to \Lien \to \Lieg \stackrel{\pi} \to \Lieq \to 0$ is said to be \emph{$\Lie$-central} if
$\Lien \subseteq Z_{\Lie}(\Lieq)$, equivalently, $[\Lien,\Lieg]_{\Lie} = 0$.

Following \cite{CK, CVDL1}, given a free presentation $0 \to \Lier \to \Lief \stackrel{\rho} \to \Lieg \to 0$ of Leibniz algebra $\Lieg$, the second \Lie-homology with trivial coefficients  is given by
 \begin{equation} \label{Hopf}
 \HL^{\Lie}_2(\Lieg) \cong \frac{\Lier \cap [\Lief, \Lief]_{\Lie}}{[\Lier, \Lief]_{\Lie}}.
  \end{equation}
  On the other hand,  the first \Lie-homology with trivial coefficients is given by  $\HL^{\Lie}_1(\Lieg) \cong \Lieg_{\Lie}$.


\section{$\Lie$-isoclinic Leibniz algebras} \label{Lie iso Lb alg}

Consider the $\Lie$-central extensions $(g) : 0 \to \Lien \stackrel{\chi} \to \Lieg \stackrel{\pi} \to \Lieq \to 0$ and $(g_i) : 0 \to \Lien_i \stackrel{\chi_i}\to \Lieg_i \stackrel{\pi_i} \to \Lieq_i \to 0, i=1, 2,$

Let be $C : \Lieq \times \Lieq \to [\Lieg, \Lieg]_{\Lie}$ given by $C(q_1,q_2)=[g_1,g_2]+[g_2,g_1]$, where $\pi(g_j)=q_j, j= 1, 2$, the $\Lie$-commutator map associated to the extension $(g)$. In a similar way is defined the $\Lie$-commutator map $C_i$ corresponding to the extensions $(g_i), i = 1, 2$.

Note that if $\Lieq$ is a Lie algebra, then $\pi([\Lieg, \Lieg]_{\Lie})=0$, hence $[\Lieg, \Lieg]_{\Lie} \subseteq \Lien \equiv \chi(\Lien)$.

\begin{De} \label{isoclinic}
The $\Lie$-central extensions $(g_1)$ and $(g_2)$ are said to be \Lie-isoclinic when there exist isomorphisms $\eta : \Lieq_1 \to \Lieq_2$ and $\xi : [\Lieg_1, \Lieg_1]_{\Lie} \to [\Lieg_2, \Lieg_2]_{\Lie}$ such that the following diagram is commutative:
\begin{equation}  \label{square isoclinic}
\xymatrix{
\Lieq_1 \times \Lieq_1 \ar[r]^{C_1} \ar[d]_{\eta \times \eta} & [\Lieg_1, \Lieg_1]_{\Lie} \ar[d]^{\xi}\\
\Lieq_2 \times \Lieq_2 \ar[r]^{C_2} & [\Lieg_2, \Lieg_2]_{\Lie}
}
\end{equation}

The pair $(\eta, \xi)$ is called a \Lie-isoclinism from $(g_1)$ to $(g_2)$ and will be denoted by $(\eta, \xi) : (g_1) \to (g_2)$.
\end{De}

\begin{Rem}
In fact it suffices to assume in Definition \ref{isoclinic} that $\xi$ is a monomorphism. Furthermore, $\eta$ determines $\xi$ uniquely.
\end{Rem}

Let $\Lieq$ be Leibniz algebra, then we can construct the following $\Lie$-central extension
\begin{equation} \label{Lie central extension}
(e_q) : 0 \to Z_{\Lie}(\Lieq) \to \Lieq \stackrel{pr_{\Lieq}} \to \Lieq/Z_{\Lie}(\Lieq) \to 0.
\end{equation}

\begin{De}
Let $\Lieg$ and $\Lieq$ be Leibniz algebras. Then $\Lieg$ and $\Lieq$ are said to be $\Lie$-isoclinic when $(e_g)$ and $(e_q)$ are $\Lie$-isoclinic $\Lie$-central extensions.

A $\Lie$-isoclinism $(\eta, \xi)$ from $(e_g)$ to $(e_q)$ is also called a $\Lie$-isoclinism  from $\Lieg$ to $\Lieq$, denoted by $(\eta, \xi) : \Lieg \sim \Lieq$.
\end{De}

\begin{Pro} \label{Lie-isoclinism}
For a $\Lie$-isoclinism $(\eta, \xi) : (g_1) \sim (g_2)$, the following statements hold:
\begin{enumerate}
\item[a)] $\eta$ induces an isomorphism $\eta' : \Lieg_1/Z_{\Lie}(\Lieg_1) \to \Lieg_2/Z_{\Lie}(\Lieg_2)$, and $(\eta', \xi)$ is a $\Lie$-isoclinism from $\Lieg_1$ to $\Lieg_2$.
    \item[b)] $\chi_1(\Lien_1) = Z_{\Lie}(\Lieg_1)$ if and only if $\chi_2(\Lien_2) = Z_{\Lie}(\Lieg_2)$.
\end{enumerate}
\end{Pro}
{\it Proof.}  Let $(\eta, \xi) : (g_1) \sim (g_2)$ be a $\Lie$-isoclinism. Then consider the map  $\phi:\Lieg_1 \to \Lieg_2/Z_{\Lie}(\Lieg_2)$ defined by $g_1\mapsto \bar{g}_2$ where $\pi_2(g_2)=(\eta\circ\pi_1)(g_1).$ Clearly, $\phi$ is a well-defined onto homomorphism since $(g_2)$ is a $\Lie$-central extension. We claim that $\Ker(\phi)=Z_{\Lie}(\Lieg_1).$ Indeed, Let $g_1\in Z_{\Lie}(\Lieg_1)$ and let $y\in\Lieg_2.$ Then set $\bar{g}_2:=\phi(g_1)$ with $g_2\in\Lieg_2.$ Now since $\phi$ is onto, $ \bar{y}=\phi(x)$ for some $x\in\Lieg_1.$
Then by commutativity of the diagram induced by   the $\Lie$-isoclinism $(\eta, \xi),$

$$\begin{aligned}
 ~[g_2,y]+[y,g_2]&= C_2(\pi_2(g_2),\pi_2(y))\\&=C_2(\eta(\pi_1(g_1)),\eta(\pi_1(x)))\\&=C_2(\eta\times\eta(\pi_1(g_1),\pi_1(x)))\\&=\xi(C_1(\pi_1(g_1),\pi_1(x)))\\&=\xi([g_1,x]+[x,g_1])\\&=\xi(0) \\&=0.\end{aligned}$$
 So $g_2\in Z_{\Lie}(\Lieg_2).$ Therefore $g_1\in \Ker(\phi).$

 Conversely, let $g_1\in \Ker(\phi)$ i.e. $\bar{g}_2:=\phi(g_1)=0$ with $g_2\in Z_{\Lie}(\Lieg_2).$ For every $x\in \Lieg_1,$  $(\eta\circ\pi_1)(x)\in \Lieq_2.$ As $\pi_2$ is onto, $(\eta\circ\pi_1)(x)=\pi_2(y)$ for some $y\in\Lieg_2.$ Again  by commutativity of the diagram induced by   the $\Lie$-isoclinism $(\eta, \xi),$
$$\begin{aligned}
 \xi(C_1(\pi_1(g_1),\pi_1(x)))&=C_2(\eta(\pi_1(g_1)),\eta((\pi_1(x)))\\&= C_2(\pi_2(g_2),\pi_2(y))\\&=[g_2,y]+[y,g_2]\\&=0.\end{aligned}$$
So $[g_1,x]+[x,g_1]=C_1(\pi_1(g_1),\pi_1(x))\in \Ker(\xi)={0}.$ Thus
 $[g_1,x]+[x,g_1]=0.$ Hence $g_1\in Z_{\Lie}(\Lieg_1).$
$\eta'$ is obtained by the first isomorphism theorem and an easy checking shows that $\eta'$  is an isomorphism.

It remains to show that the following diagram is commutative:
 \[
 \xymatrix{
 \Lieg_1/Z_{\Lie}(\Lieg_1) \times \Lieg_1/Z_{\Lie}(\Lieg_1) \ar[r]^{\quad \quad \quad \bar{C}_1} \ar[d]_{\eta' \times \eta'} & [\Lieg_1, \Lieg_1]_{\Lie} \ar[d]^{\xi}\\
 \Lieg_2/Z_{\Lie}(\Lieg_2) \times \Lieg_2/Z_{\Lie}(\Lieg_2) \ar[r]^{\quad \quad \quad \bar{C}_2} & [\Lieg_2, \Lieg_2]_{\Lie}
 }
 \]
 Indeed, since $\Lien_i \subseteq Z_{\Lie}(\Lieq_i)$ for $i=1,2,$ we have  well-defined maps $P_i: \Lieq_i\to \Lieg_i/Z_{\Lie}(\Lieg_i),~z^i\mapsto \bar{x}^i$ where $z^i=\pi_i(x^i).$  These maps induce the well-defined maps $\bar{C}_i:\Lieg_i/Z_{\Lie}(\Lieg_i) \times \Lieg_i/Z_{\Lie}(\Lieg_i) \to [\Lieg_i, \Lieg_i]_{\Lie}$ defined by $\bar{C}_i(\bar{x}^i_1,\bar{x}^i_2):=C_1(z_1^i,z_2^i).$  Now let  $g_1,g_2\in\Lieg_1$ such that $\pi_1(g_1)=q_1$ and $\pi_2(g_2)=q_2$, and let $h_1, h_2 \in \Lieg_2$ such that $\pi_2(h_1) = (\eta \circ \pi_1)(g_1)$ and $\pi_2(h_2) = (\eta \circ \pi_1)(g_2)$. Then
  $$
  \begin{aligned}
\xi\circ\bar{C}_1(\bar{g}_1,\bar{g}_2)&=\xi\circ C_1(q_1,q_2)=(C_2\circ(\eta\times\eta)(q_1,q_2)))\\&=C_2((\eta\circ\pi_1)(g_1),(\eta\circ\pi_1)(g_2))\\ & = C_2(\pi_2(h_1), \pi_2(h_2)) \\ & = \bar{C}_2(\bar{h}_1, \bar{h}_2) \\ &=\bar{C}_2(\eta'(\bar{g}_1),\eta'(\bar{g}_2))\\&=(\bar{C}_2\circ(\eta'\times\eta'))(\bar{g}_1,\bar{g}_2).
 \end{aligned}
 $$  Therefore the diagram above is commutative.

 To show b), assume that $\chi_1(\Lien_1) = Z_{\Lie}(\Lieg_1)$ and let $z\in\chi_2(\Lien_2) $ i.e. $\pi_2(z)=0.$ Since $\eta'$ is onto, $\eta'(\bar{z}')=\bar{z}$ for some $z'\in\Lieg_1.$ Then $(\eta\circ\pi_1)(z')=\pi_2(z)=0$ i.e. $\eta(\pi_1(z'))=0.$ This implies that $\pi_1(z')=0$ since $\eta$ is one to one. So $z'\in\chi_1(\Lien_1) = Z_{\Lie}(\Lieg_1),$ and thus $\bar{z}=\eta'(\bar{z}')=0$ in $\Lieg_2/Z_{\Lie}(\Lieg_2).$ Hence $z\in Z_{\Lie}(\Lieg_2).$  Conversely, let $z\in Z_{\Lie}(\Lieg_2).$ Then $\bar{z}=0$ in  $\Lieg_2/Z_{\Lie}(\Lieg_2).$ Since $\eta'$ is an isomorphism, $\eta'(\bar{t})=0$ for some $t\in Z_{\Lie}(\Lieg_1)=\chi_1(\Lien_1).$ So $\pi_2(z)=(\eta\circ\pi_1)(t)=\eta(\pi_1(t))=\eta(0)=0.$ Therefore $z\in\chi_2(\Lien_2).$ \\
 The other implication is identical. \rdg

\begin{Pro} \label{equivalence relation}
$\Lie$-isoclinism is an equivalence relation.
\end{Pro}
{\it Proof.} Direct cheking. \rdg
\bigskip

$\Lie$-isoclinisms of the form $(\eta, \xi) : (g) \sim (g)$, respectively $(\eta, \xi) : \Lieg \sim \Lieg$, are called \emph{autoclinisms} of $(g)$, respectively $\Lieg$, and they constitute the group $Autcl(g)$, respectively $Autcl(\Lieg)$. Each $\Lie$-isoclinism from $(g_1)$ to $(g_2)$ induces an isomorphism from $Autcl(g_1)$ to $Autcl(g_2)$.

Let $(\alpha, \beta, \gamma): (g_1) \to (g_2)$ be an isomorphism of $\Lie$-central extensions (i.e. $\alpha, \beta, \gamma$ are bijective Leibniz homomorphisms such that $\chi_2 \circ \alpha= \beta \circ \chi_1$ and $\pi_2 \circ \beta = \gamma \circ \pi_1$) and $\beta'$ the restriction of $\beta$ to $[\Lieg_1, \Lieg_1]_{\Lie}$. Then $(\gamma, \beta')$ is a $\Lie$-isoclinism from $(g_1)$ to $(g_2)$.

From Proposition \ref{Lie-isoclinism}, the classes of $\Lie$-isoclinic Leibniz algebras can be regarded as those isoclinism classes that consist of $\Lie$-central extensions $(g)$ with $\chi(\Lien)=Z_{\Lie}(\Lieg)$, i. e. $(g)$ is isomorphic to $(e_g)$.

The Leibniz algebras which are $\Lie$-isoclinic to the trivial Leibniz algebra are exactly the abelian ones. Indeed, if
$\Lieg_1$ is a Leibniz algebra and $\Lieg_2$ the trivial Leibniz algebra which are $\Lie$-isoclinic, then  the  isomorphism $\xi : [\Lieg_1, \Lieg_1]_{\Lie} \to [\Lieg_2, \Lieg_2]_{\Lie}$ and the fact $[\Lieg_2, \Lieg_2]_{\Lie} =0$ implies that $[\Lieg_1, \Lieg_1]_{\Lie}=0$, so the isomorphism  $\eta : \Lieg_1 \to \Lieg_2$  implies that $\Lieg_1$ is an abelian Leibniz algebra.

Conversely,  every abelian Leibniz algebra $\Lieg$ is {\Lie}-isoclinic to the trivial Leibniz algebra on the underlying set $\Lieg.$ Indeed, for both Leibniz algebras, $Z_{\Lie}(\Lieg)=\Lieg$ and thus $\Lieg/Z_{\Lie}(\Lieg)=0.$ We then have the {\Lie}- isoclinism $(\eta=id_{\Lieg},\xi=0):(\Lieg)\sim(\Lieg).$

\begin{Ex}
Following we provide an example of two non-abelian \Lie-isoclinic Leibniz algebras.

Consider the non-abelian complex Leibniz algebras   $\Lieg_1 = span\{e_1, e_2\}$  with non-zero multiplications $[e_1, e_1]=[e_2, e_1] = e_2$ (it belongs to the third isomorphism class of Lemma 2 in \cite{Cu}),  and   $\Lieg_2 = span\{a_1, a_2, a_3\}$ with non-zero multiplications $[a_1, a_1]=[a_2, a_1]=[a_3,a_1]= a_3$ (it belongs to the isomorphism class 2 (d) of the classification given  in \cite{CILL}). Clearly $[\Lieg_1,\Lieg_1]_{\Lie}=span\{e_2\},$ $[\Lieg_2,\Lieg_2]_{\Lie}=span\{a_3\},$ $Z_{\Lie}(\Lieg_1)=0$ and $Z_{\Lie}(\Lieg_2)=span\{a_2-a_3\}.$ So $\Lieq_1:=\Lieg_1/Z_{\Lie}(\Lieg_1)\simeq \Lieg_1.$ Also applying the first isomorphism theorem on the map $\phi: \Lieg_2\to\Lieg_2$ defined by $\phi(a_1)=a_1$ and $\phi(a_2)=\phi(a_3)=a_3$ shows $\Lieq_2:=\Lieg_2/Z_{\Lie}(\Lieg_2)\simeq span\{a_1, a_3\}.$ Now define $\eta : \Lieq_1 \to \Lieq_2$ by $\eta(e_1)=a_1$ and $\eta(e_2)=a_3,$ and $\xi : [\Lieg_1, \Lieg_1]_{\Lie} \to [\Lieg_2, \Lieg_2]_{\Lie}$ by $\xi(e_2)=a_3.$ It is easy to show that $(\eta,\xi):(\Lieg)\sim(\Lieg)$ is a {\Lie}-isoclinism.
\end{Ex}

\begin{Pro} \label{backward}
Let $(\eta, \xi) : (g_1) \sim (g_2)$ be a $\Lie$-isoclinism. Then the backward induced extension $(\eta^{\ast}(g_2)): 0 \to \Lien_2 \to \Lieg_2^{\eta} \stackrel{\overline{\pi}_2}\to \Lieq_1 \to 0$, obtained by pulling back along $\eta$ (where $\Lieg_2^{\eta} = \{ (g,q) \in  \Lieg_2 \times  \Lieq_1 \mid \pi_2(g) = \eta(q) \}$) is a $\Lie$-central extension isomorphic to $(g_2)$, and $(id_{\Lieq_1}, \xi) : (g_1) \sim (\eta^{\ast}(g_2))$ is a $\Lie$-isoclinism.
\end{Pro}
{\it Proof.}
The first statement only requires a routine checking.

To show that $(id_{\Lieq_1}, \xi) : (g_1) \sim (\eta^{\ast} g_2)$ is a $\Lie$-isoclinism, have in mind that  $\xi$ induces an isomorphism (call it also $\xi$) $\xi:[\Lieg_1, \Lieg_1]_{\Lie}\to[\Lieg_2^{\eta}, \Lieg_2^{\eta}]_{\Lie}$
 defined by $\xi([g_1,g_1']+[g_1',g_1])=([g_2,g_2']+[g_2',g_2],[q_1,q_1']+[q_1',q_1])$, where  $\pi_2(g_2)=\eta(q_1) = \eta \pi_1(g_1)$. The commutativity of the diagram
 \[
 \xymatrix{
 \Lieq_1 \times \Lieq_1 \ar[r]^{C_1} \ar[d]_{id_{\Lieq_1} \times id_{\Lieq_1}} & [\Lieg_1, \Lieg_1]_{\Lie} \ar[d]^{\xi}\\
 \Lieq_1 \times \Lieq_1 \ar[r]^{C_2^{\eta}} & [\Lieg_2^{\eta}, \Lieg_2^{\eta}]_{\Lie}
 }
 \]
 where $C_2^{\eta}(q_1,q_1')=([g_2,g_2']+[g_2',g_2],[q_1,q_1']+[q_1',q_1])$, follows directly.
 \rdg
\bigskip

Thanks to Proposition \ref{backward}, in many cases we can restrict ourselves to $\Lie$-isoclinisms $(\eta, \xi)$ in which $\eta$ is an identity map.

\begin{Pro}
Let $(\eta, \xi) : (g_1)  \sim (g_2)$ be an $\Lie$-isoclinism. Then the following statements hold:
\begin{enumerate}
\item[a)] $\pi_2 \circ  \xi(g) = \eta \circ \pi_1(g)$, for all $g \in [\Lieg_1, \Lieg_1]_{\Lie}$.
\item[b)] $\xi \left( \chi_1(\Lien_1) \cap [\Lieg_1, \Lieg_1]_{\Lie} \right) = \chi_2(\Lien_2) \cap [\Lieg_2, \Lieg_2]_{\Lie}$.
\item[c)] $\xi([g,x]+[x,g]) = [h,\xi(x)]+[\xi(x),h]$, for all $x \in  [\Lieg_1, \Lieg_1]_{\Lie}, g \in \Lieg_1, h \in \Lieg_2$ such that $\eta \circ \pi_1(g) = \pi_2(h)$.
\end{enumerate}
\end{Pro}
{\it Proof.}
{\it a)}  Proposition \ref{backward} establishes that the  $\Lie$-isoclinism $(\eta, \xi) : (g_1)  \sim (g_2)$ induces a $\Lie$-isoclinism $(id_{\Lieq_1}, \xi) : (g_1) \sim (\eta^{\ast} g_2) $ where $(\xi(g),\pi_1(g))\in [\Lieg_2^{\eta}, \Lieg_2^{\eta}]_{\Lie}$  for all $g\in[\Lieg_1, \Lieg_1]_{\Lie} \subseteq \Lieg_2^{\eta}$.

{\it b)} Let $g\in\chi_1(\Lien_1)\cap [\Lieg_1, \Lieg_1]_{\Lie}$,   then $\pi_1(g)=0$ which implies by {\it a)} that $\pi_2(\xi(g))=0.$ So $\xi(g)\in \Ker(\pi_2)=\chi_2(\Lien_2).$

Conversely, let $h\in \chi_2(\Lien_2) \cap [\Lieg_2, \Lieg_2]_{\Lie}.$ Then $\pi_2(h)=0$ and  $h=\xi(g)$ for some $g\in [\Lieg_1, \Lieg_1]_{\Lie}$ since $\xi$ is onto. Again by {\it a)} we have $\eta \circ \pi_1(g)=0.$ This implies that $\pi_1(g)=0$.  Thus $g\in\chi_1(\Lien_1).$

 {\it c)} Let $x \in  [\Lieg_1, \Lieg_1]_{\Lie}, g \in \Lieg_1$ and $h \in \Lieg_2$ such that $\eta (\pi_1(g)) = \pi_2(h)$.
Then $$
\begin{aligned}
\xi([g,x]+[x,g]) &= \xi(C_1(\pi_1(g),\pi_1(x)))\\&=C_2(\eta \circ \pi_1(g),\eta \circ \pi_1(x))\\&=C_2(\pi_2(h),\pi_2(\xi(x))\\&=[h,\xi(x)]+[\xi(x),h].
\end{aligned}
$$ \rdg
\bigskip

As we have mentioned above, isomorphisms of $\Lie$-central extensions induce $\Lie$-isoclinisms. This observation gives rise to the following:

\begin{De} \label{isoclinic1}
A homomorphism of  $\Lie$-central extensions $(\alpha, \beta, \gamma) : (g_1) \to (g_2)$ is said to be $\Lie$-isoclinic, if there exists an isomorphism $\beta' : [\Lieg_1,\Lieg_1]_{\Lie} \to [\Lieg_2,\Lieg_2]_{\Lie}$ with $(\gamma, \beta') : (g_1) \sim (g_2)$.

If $\beta$ is in addition an epimorphism (resp., monomorphism), then  $(\alpha, \beta, \gamma)$ is called an isoclinic epimorphism (resp., monomorphism).
\end{De}

\begin{Pro} \label{equivalence}
For a homomorphism of $\Lie$-central extensions $(\alpha, \beta, \gamma) : (g_1) \to (g_2)$, the following statements hold:
\begin{enumerate}
\item[a)]  $(\alpha, \beta, \gamma)$ is  $\Lie$-isoclinic if and only if $\gamma$ is an isomorphism and $\Ker(\beta) \cap [\Lieg_1,\Lieg_1]_{\Lie} =0$.
    \item[b)] If  $(\alpha, \beta, \gamma)$ is  $\Lie$-isoclinic and $\beta'$ as in Definition \ref{isoclinic1}, then $\beta' = \beta_{\mid  [\Lieg_1,\Lieg_1]_{\Lie}}$.
\end{enumerate}
\end{Pro}
{\it Proof.} {\it a)} Assume that $(\alpha, \beta, \gamma) : (\Lieg_1) \to (\Lieg_2)$ is $\Lie$-isoclinic, then $(\gamma, \beta') : (\Lieg_1) \sim (\Lieg_2)$ is a  $\Lie$-isoclinism for some isomorphism $\beta' : [\Lieg_1,\Lieg_1]_{\Lie} \to [\Lieg_2,\Lieg_2]_{\Lie}.$ This implies by definition that $\gamma$ is an isomorphism. Now let $m\in \Ker(\beta) \cap [\Lieg_1,\Lieg_1]_{\Lie}.$ Then $\beta(m)=0$ and $m=[g_1,g_2]+[g_2,g_1]$ for some $g_1,g_2\in \Lieg_1.$ Also let $q_1=\pi_1(g_1)$ and $q_2=\pi_1(g_2).$
Since $(\gamma, \beta') : (\Lieg_1) \sim (\Lieg_2)$ is a  $\Lie$-isoclinism, we have
$$
\begin{aligned}
\beta'(m)&=\beta'([g_1,g_2]+[g_2,g_1])\\&=\beta'(C_1(\pi_1(g_1),\pi_1(g_2)))\\&=(\beta'\circ C_1)(\pi_1(g_1),\pi_1(g_2))\\&=(C_2\circ(\gamma\times\gamma))((\pi_1(g_1),\pi_1(g_2))\\&=C_2(\gamma(\pi_1(g_1)),\gamma(\pi_1(g_2)))\\&=C_2(\pi_2(\beta(g_1)),\pi_2(\beta(g_2))) \\&=[\beta(g_1),\beta(g_2)]+[\beta(g_2),\beta(g_1)]\\&=\beta([g_1,g_2]+[g_2,g_1])\\&=\beta(m)=0.
\end{aligned}$$
Since $\beta'$ is one-to-one, it follows that $m=0.$

Conversely, assume that $\Ker(\beta) \cap [\Lieg_1,\Lieg_1]_{\Lie} =0.$ Define $\beta':[\Lieg_1,\Lieg_1]_{\Lie}\to [\Lieg_2,\Lieg_2]_{\Lie}$ by $\beta'(g)=\beta(g)$, which is  one-to-one. It remains to show that $\beta'$ is onto. Let $y\in [\Lieg_2,\Lieg_2]_{\Lie}.$ Then $y=[g_2,g_2']+[g_2',g_2]$ for some $g_2, g'_2\in\Lieg_2.$ Since $\pi_1$ and $\gamma_1$ are onto, it follows that $\pi_2(g_2)=(\gamma\circ\pi_1)(g_1)$ and $\pi_2(g'_2)=(\gamma\circ\pi_1)(g'_1)$ for some $g_1,g_1'\in\Lieg_1.$ By the homomorphism $(\alpha,\beta,\gamma)$, we have      $(\gamma\circ\pi_1)(g_1)=(\pi_2\circ\beta)(g_1)$ and $(\gamma\circ\pi_1)(g_1')=(\pi_2\circ\beta)(g_1'),$ which implies that  $g_2-\beta(g_1)=\chi_2(n)$ and $g_2'-\beta(g_1')=\chi_2(n')$ for some $n,n'\in\Lien_2.$  We now have
$$
\begin{aligned}
\beta'([g_1,g_1']+[g_1',g_1])&=\beta([g_1,g_1']+[g_1',g_1])\\&=[\beta(g_1),\beta(g_1')]+[\beta(g_1'),\beta(g_1)]\\
&=[g_2-\chi_2(n),g_2'-\chi_2(n')]+[g'_2-\chi_2(n'),g_2-\chi_2(n)]\\&=y-([\chi_2(n),g_2'-\chi_2(n')]+[g'_2-\chi_2(n'),\chi_2(n)]) \\
&~~~~~~-([g_2,\chi_2(n')]+[\chi_2(n'),g_2'])\\&=y.
\end{aligned}$$

{\it b)} follows directly from the proof of {\it a)}.
 \rdg

\begin{Pro} \label{character}
Let $\beta : \Lieg \to \Lieh$ be a homomorphism of Leibniz algebras. Then $\beta$ induces a $\Lie$-isoclinic homomorphism from $(e_{\Lieg})$ to $(e_{\Lieh})$ if and only if $\Ker(\beta) \cap [\Lieg,\Lieg]_{\Lie}=0$ and $\Image(\beta) + Z_{\Lie}(\Lieh) = \Lieh$.

In this case we call $\beta$ a $\Lie$-isoclinic homomorphism.
\end{Pro}
{\it Proof.} Assume that $\Ker(\beta) \cap [\Lieg,\Lieg]_{\Lie}=0$ and $\Image(\beta)+ Z_{\Lie}(\Lieh) = \Lieh.$ First we prove that  $\beta(Z_{\Lie}(\Lieg))\subseteq Z_{\Lie}(\Lieh).$ Indeed let $h\in\Lieh.$ Since $\Image(\beta)+ Z_{\Lie}(\Lieh) = \Lieh,$ it follows that $h=\beta(x)+h_0$ for some $x\in\Lieg$ and $h_0\in Z_{\Lie}(\Lieh).$ Then
$$
\begin{aligned}
~~[\beta(g),h]+[h,\beta(g)]&= [\beta(g),\beta(x)+h_0]+[\beta(x)+h_0,\beta(g)]\\&=[\beta(g),\beta(x)]+[\beta(x),\beta(g)]\\&=0.
\end{aligned}
$$
So the maps $\alpha:=\beta|_{Z_{\Lie}(\Lieg)}$  and $\gamma:\Lieg/Z_{\Lie}(\Lieg)\to \Lieh/Z_{\Lie}(\Lieh)$ given by $\gamma(\overline{{g}})=\overline{{\beta(g)}}$ are well-defined homomorphisms and it is readily verified that $(\alpha,\beta,\gamma): (e_{\Lieg})\to (e_{\Lieh})$  is a homomorphism of \Lie-central extensions. To show that it is ${\Lie}$-isoclinic, it is enough to show by Proposition \ref{equivalence} that $\gamma$  is an isomorphism. To  show that $\gamma$ is one-to-one, let $g\in\Lieg$ such that $\gamma(\overline{g})=0.$ Then $\beta(g)\in Z_{\Lie}(\Lieh).$ We claim that $g\in Z_{\Lie}(\Lieg).$ Indeed if $g\notin Z_{\Lie}(\Lieg),$ then $m:=[g,g']+[g',g]\neq 0$ for some $g'\in\Lieg.$ But $\beta(m)=[\beta(g),\beta(g')]+[\beta(g'),\beta(g)]=0$ because $\beta(g)\in Z_{\Lie}(\Lieh).$ This implies that $m\in \Ker(\beta)\cap [\Lieg,\Lieg]_{\Lie},$ and thus $m=0.$ A contradiction. Next we show that $\gamma$ is onto. Let $h\in\Lieh.$  Since $\Image(\beta)+ Z_{\Lie}(\Lieh) = \Lieh,$ it follows that $h=\beta(x)+h_0$ for some $x\in\Lieg$ and $h_0\in Z_{\Lie}(\Lieh).$ Clearly, $\overline{h}=\overline{\beta(x)}=\gamma(\overline{x}).$

Conversely, assume that $\beta$ induces a $\Lie$-isoclinic homomorphism  $(\alpha,\beta,\gamma)$ from $(e_{\Lieg})$ to $(e_{\Lieh}).$ Then again by Proposition \ref{equivalence},  $\Ker(\beta) \cap [\Lieg,\Lieg]_{\Lie}=0.$ It remains to show that $\Image(\beta)+ Z_{\Lie}(\Lieh) = \Lieh.$ Clearly $\Image(\beta)+ Z_{\Lie}(\Lieh) \subseteq \Lieh.$ Now let $h\in\Lieh.$ Following the notation in equation (\ref{Lie central extension}), $pr_{\Lieh}$ and $\gamma$ are onto,  then $pr_{\Lieh}(h)=\gamma\circ pr_{\Lieg}(g)$ for some $g\in \Lieg.$ On the other hand we have by  the homomorphism $(\alpha,\beta,\gamma)$ that $(\gamma\circ pr_{\Lieg})(g)=(pr_{\Lieh}\circ\beta)(h),$  which implies that $h-\beta(g)\in \Ker(pr_{\Lieh})=Z_{\Lie}(\Lieh).$ Therefore $h=\beta(g)+n$  for some $n\in Z_{\Lie}(\Lieh).$ This completes the proof.    \rdg

\begin{Le} \label{pull-back}
Let $P=\{(a,c) \in A \times C : \alpha(a) = \beta(c)\}$ be the pull-back of a homomorphism $\alpha: A \to B$ and an isomorphism $\beta : C \to B$ in the category of Leibniz algebras.
 \[ \xymatrix{
 P \ar[r]^{\overline{\alpha}} \ar[d]^{\overline{\beta}} & C \ar[d]^{\beta}\\
 A \ar[r]^{\alpha} & B
 }\] Then
 \begin{enumerate}
 \item[a)] $\beta$ monomorphism implies $\overline{\beta}$ monomorphism.

 \item[b)] Given $\pi : G \to C, \gamma : G \to A$ such that $\alpha \circ \gamma = \beta \circ \pi$, then there exists a unique homomorphism $\omega:G \to P$ satisfying $\overline{\alpha}\circ\omega=\pi$ and $\overline{\beta}\circ\omega=\gamma$. If $\pi$ is an epimorphism, then $\omega$ is an epimorphism.
 \end{enumerate}
 {\it Proof.} Direct checking. \rdg

\end{Le}
\begin{Pro}\
\begin{enumerate}
\item[a)] A homomorphism  $(\alpha, \beta, \gamma): (g_1) \to (g_2)$ is $\Lie$-isoclinic if and only if $\gamma$ is an isomorphism and $\beta : \Lieg_1 \to \Lieg_2$ is a $\Lie$-isoclinic homomorphism of Leibniz algebras.

    \item[b)] The composition of $\Lie$-isoclinic homomorphisms is an isoclinic homomorphism.

    \item[c)] Each $\Lie$-isoclinic homomorphism is a composition of a $\Lie$-isoclinic epimorphism and a  $\Lie$-isoclinic monomorphism.
\end{enumerate}
\end{Pro}
{\it Proof.}
{\it a)} Assume that $(\alpha, \beta, \gamma): (g_1) \to (g_2)$ is $\Lie$-isoclinic. Then by Proposition 3.10, $\gamma$ is an isomorphism and $\Ker(\beta) \cap [\Lieg_1,\Lieg_1]_{\Lie}=0.$ It remains to show that $\Image(\beta)+ Z_{\Lie}(\Lieg_2)=\Lieg_2.$ Clearly $\Image(\beta)+ Z_{\Lie}(\Lieg_2)\subseteq\Lieg_2.$ Now let $h\in\Lieg_2.$ Since $\pi_2$ and $\gamma$ are onto, it follows that  $\pi_2(h)=\gamma\circ\pi_1(g)$ for some $g\in \Lieg_1.$ Also by  the homomorphism $(\alpha,\beta,\gamma),$ we have $(\gamma\circ\pi_1)(g)=(\pi_2\circ\beta)(g).$  So $h-\beta(g)\in \Ker(\pi_2)=\Image(\chi_2).$ Therefore $h=\beta(g)+\chi_2(n)$  for some $n\in\Lien_2.$ This proves the assertion since $\chi_2(n)\in Z_{\Lie}(\Lieg_2)$ as $[\Lien_2,\Lieg_2]_{\Lie}=0.$

Conversely, assume that $\beta$ is a {\Lie}-isoclinic homomorphism of Leibniz algebras. Then $\Ker(\beta) \cap [\Lieg_1,\Lieg_1]_{\Lie}=0$ and $\Image(\beta)+ Z_{\Lie}(\Lieg_2)=\Lieg_2.$ Again by Proposition 3.10, it is enough to show that $\gamma$ is an isomorphism.
To show that $\gamma$ is onto, let $q\in\Lieq_2$ and let $h\in\Lieg_2$ such that $q=\pi_2(h).$ Then  $h=\beta(x)+h_0$ for some $x\in\Lieg_1$ and $h_0\in Z_{\Lie}(\Lieg_2).$ We then have $\gamma(\pi_1(x))=$
$\pi_2(\beta(x))=\pi_2(\beta(x)+h_0)=\pi_2(h)=q.$ To show that $\gamma$ is one-to-one, let $q_1\in \Lieq_1$ with $\gamma(q_1)=0,$ i.e. $\gamma(\pi_1(g_1))=0$ for some $g_1\in\Lieg_1.$   By the homomorphism $(\alpha,\beta,\gamma),$ it follows that $\pi_2(\beta(g_1))=0$ i.e. $\beta(g_1)=\chi_2(n)$ for some $n\in\Lien_2.$ We claim that $g_1\in \Lien_1$ in which case  $q_1=\pi_1(g_1)=0.$ Otherwise,  there exists  $g_1'\in\Lieg_1$ such that $m:=[g_1,g_1']+[g_1',g_1]\neq0.$ So
$$
\begin{aligned}
\beta(m)&=[\beta(g_1),\beta(g_1')]+[\beta(g_1'),\beta(g_1)]\\&=
[\chi_2(n),\beta(g_1')]+[\beta(g'_1),\chi_2(n)]\\&=0.
\end{aligned}$$
This contradicts the fact that $\Ker(\beta) \cap [\Lieg_1,\Lieg_1]_{\Lie}=0.$

The proof of {\it b)} is straightforward.

{\it c)} Let $\beta: \Lieg \to \Lieh$ be a \Lie-isoclinic homomorphism. Then $\beta$ induces a \Lie-isoclinism $(\alpha, \beta, \gamma) : (e_{\Lieg}) \to (e_{\Lieh})$ where $\alpha=\beta_{\mid Z_{\Lie}(\Lieg)}, \gamma = \overline{\beta}$, that is $\gamma(\overline{g})=\overline{\beta(g)}$. Moreover, by Proposition 3.10, we know that $\gamma$  is an isomorphism and $\Ker(\beta) \cap [\Lieg, \Lieg]_{\Lie}=0$. Equivalently, by Definition 3.9, $\beta' :   [\Lieg, \Lieg]_{\Lie} \to [\Lieh, \Lieh]_{\Lie}$ is an isomorphism and $(\gamma, \beta'): (e_{\Lieg}) \sim (e_{\Lieh})$. Now applying Proposition 3.7 we obtain the decomposition
$$(e_{\Lieg}) \stackrel{(\overline{\omega}_{\mid}, \omega, 1)}\to (e_{\Lieh^{\gamma}}) \stackrel{(1, \overline{\gamma}, \gamma)}\to (e_{\Lieh})$$
where $w$ is induced by pull-back properties.

 In fact we have the following diagram:
\[ \xymatrix{
 e_{g} :  0 \ar[r] & Z_{\Lie}(\Lieg) \ar[r]^{} \ar[d] & \Lieg \ar[r]^{\pi_1 \quad \quad} \ar[d]^{\omega} & \Lieg/Z_{\Lie}(\Lieg)  \ar@{=}[d] \ar[r]& 0 \\
(e_{h^{\gamma}}) :  0 \ar[r] & Z_{\Lie}(\Lieh) \ar@{=}[d] \ar[r] & \Lieh^{\gamma}\ar[r]^{\overline{{\pi}}_2 \quad \quad} \ar[d]^{\overline{\gamma}} & \Lieg/Z_{\Lie}(\Lieg) \ar[r] \ar[d]^{\gamma}& 0\\
e_{h} :  0 \ar[r] & Z_{\Lie}(\Lieh) \ar[r]^{}  & \Lieh \ar[r]^{\pi_2 \quad \quad}  & \Lieh/Z_{\Lie}(\Lieh)\ar[r] & 0
 }\]
Now Lemma \ref{pull-back} implies that $\omega$ is an epimorphism and $\overline{\gamma}$ is a monomorphism, as required. \rdg

\bigskip

In the following consider an isomorphism $\eta : \Lieq_1 \to \Lieq_2$ and put $\Lieq = \Lieq_1$. Construct the following backward (see Proposition \ref{backward}) and product $\Lie$-central extensions:
\[
\xymatrix{
(\eta^{\ast} g_2) :  0 \ar[r] & \Lien_2 \ar[r]^{\chi_2} \ar@{=}[d] & \Lieg_2^{\eta} \ar[r]^{\overline{\pi}_2 \quad } \ar[d]^{\overline{\eta}} & \Lieq_1 = \Lieq \ar[d]^{\eta} \ar[r]& 0 \\
(g_2) :  0 \ar[r] & \Lien_2 \ar[r]^{\chi_2}  & \Lieg_2 \ar[r]^{\pi_2}  & \Lieq_2 \ar[r] & 0
}
\]

\[
\xymatrix{
(\widetilde{g}) :  0 \ar[r] & \Lien_1 \times \Lien_2 \ar[r]^{\lambda} \ar@{=}[d] & \widetilde{\Lieg} \ar[r]^{\rho} \ar[d]^{\overline{\Delta}} &  \Lieq \ar[d]^{\Delta} \ar[r]& 0 \\
(g_1 \times \eta^{\ast} g_2) :  0 \ar[r] & \Lien_1 \times \Lien_2 \ar[r]^{\chi_1 \times \chi_2}  & \Lieg_1 \times \Lieg_2^{\eta} \ar[r]^{\quad \pi_1 \times \overline{\pi}_2}  & \Lieq \times \Lieq \ar[r] & 0
}
\]
where $\Delta(q) = (q, q), q \in \Lieq$ is the diagonal map and $\widetilde{\Lieg} =(\Lieg_1 \times \Lieg_2^{\eta})^{\Delta}$.

An easy computation shows that $ \widetilde{\Lieg} \equiv \{ (g_1, g_2), g_i \in \Lieg_i, i=1,2 \mid \eta \circ \pi_1(g_1) = \pi_2(g_2) \}, \rho(g_1, g_2) = \pi_1(g_1), \lambda \equiv \chi_1 \times \chi_2$ and $(\widetilde{g})$ is a $\Lie$-central extension.

Furthermore we denote by $\sigma_i : \Lien_1 \times \Lien_2 \to \Lien_i$ and $\tau_i : \widetilde{\Lieg} \to \Lieg_i$ the $i$-th projection ($i=1, 2$) and put $\gamma_1 = id_{\Lieq} : \Lieq \to \Lieq_1$, $\gamma_2 = \eta : \Lieq \to \Lieq_2$.

\begin{Pro} \label{Lie iso epi}\
\begin{enumerate}
\item[a)] The following diagrams are commutative for $i=1, 2$:
\begin{equation} \label{tilda g}
\xymatrix{
(\widetilde{g}) :  0 \ar[r] & \Lien_1 \times \Lien_2 \ar[r]^{\quad \lambda} \ar[d]^{\sigma_i} & \widetilde{\Lieg} \ar[r]^{\rho} \ar[d]^{\tau_i} &  \Lieq \ar@{>->}[d]^{\gamma_i} \ar[r]& 0 \\
(g_i) :  0 \ar[r] & \Lien_i \ar[r]^{\chi_i}  & \Lieg_i \ar[r]^{\pi_i}  & \Lieq_i \ar[r] & 0
}
\end{equation}
\item[b)] The isomorphism $\eta : \Lieq_1 \to \Lieq_2$ induces a $\Lie$-isoclinism from $(g_1)$ to $(g_2)$ if and only if $(\sigma_i, \tau_i, \gamma_i), i = 1, 2$, are $\Lie$-isoclinic epimorphisms.
\end{enumerate}
\end{Pro}
{\it Proof.} {\it a)} Let $(g_1,g_2)\in\widetilde{\Lieg}.$ Then for $i=1,2,$ we have
 $$(\gamma_i\circ\rho)(g_1,g_2) = \gamma_i(\pi_1(g_1))=\begin{cases} \pi_1(g_1) ,& \mbox{if} ~~i=1\mbox{}\\\eta(\pi_1(g_1)),& \mbox{if}~~i=2 \mbox{}\end{cases} = (\pi_i\circ\tau_i)(g_1,g_2).$$
 On the other hand,  let $(n_1,n_2)\in \Lien_1\times\Lien_2.$ Then for $i=1,2,$
 $$(\chi_i\circ\sigma_i)(n_1,n_2)=\chi_i(n_i)=\tau_i(\chi(n_1),\chi(n_2))= (\tau_i\circ\lambda)(n_1,n_2).$$

{\it b)} Assume that $(\sigma_i,\tau_i,\gamma_i),~i=1,2$ are $\Lie$-isoclinic epimorphisms. Then by Proposition 3.10, $\Ker(\tau_i)\cap[\widetilde{\Lieg},\widetilde{\Lieg}]_{\Lie}=0,$ and $\gamma_i$  are isomorphisms, $i = 1, 2$.

Define $\eta:= \gamma_2$ and $\xi:[{\Lieg_1},{\Lieg_1}]_{\Lie}\to[{\Lieg_2},{\Lieg_2}]_{\Lie}$ by $\xi(g)=h$ where $\pi_2(h)=\eta\circ\pi_1(g).$ $\xi$ is a well-defined homomorphism and the diagram
\[
\xymatrix{
     \Lieq_1 \times \Lieq_1 \ar[r]^{C_1} \ar[d]_{\eta \times \eta} & [\Lieg_1, \Lieg_1]_{\Lie} \ar[d]^{\xi}\\
     \Lieq_2 \times \Lieq_2 \ar[r]^{C_2} & [\Lieg_2, \Lieg_2]_{\Lie}
     }
     \]
     is commutative. It remains to show that $\xi$ is bijective.

$\xi$ is one-to-one, since for any $g\in[\Lieg_1,\Lieg_1]_{\Lie}$ such that $\xi(g)=0$,  we have that   $(g,\xi(g))\in \Ker(\tau_2)\cap[\widetilde{\Lieg},\widetilde{\Lieg}]_{\Lie}=0,$ which implies that  $g=0.$

 $\xi$ is onto, since for any  $h=[b_1,b_2]+[b_2,b_1] \in[\Lieg_2,\Lieg_2]_{\Lie}$ we have that  there exists $a_j\in\Lieg_1,~j=1,2$ with $\pi_2(b_j)=\eta\circ\pi_1(a_j)$ such that  $\tau_2(a_j,b_j)=b_j$. Let be $g:=[a_1,a_2]+[a_2,a_1].$ Clearly, $\pi_2(h)=\eta\circ\pi_1(g).$ So $\pi_2(\xi(g))=\pi_2(h)$ by definition of $\xi$.  Therefore $h=\xi(g)+\chi_2(n)$ for some $n\in\Lien_2.$ However, it is obvious that $(0,\chi_2(n))\in \Ker(\tau_1)\cap[\widetilde{\Lieg},\widetilde{\Lieg}]_{\Lie}=0.$ This implies that $\chi_2(n)=0,$ and thus $h=\xi(g).$

  Conversely,
  Suppose that $\eta:\Lieq_1\to\Lieq_2$ induces a $\Lie$-isoclinism $(\eta,\xi): (g_1)\to(g_2).$ Then we have an isomorphism $\xi : [\Lieg_1, \Lieg_1]_{\Lie} \to [\Lieg_2, \Lieg_2]_{\Lie}$ defined by $\xi([g_1,g_2]+[g_2,g_1])=[h_1,h_2]+[h_2, h_1], g_i \in \Lieg_1, h_i \in \Lieg_2, \eta \circ \pi_1(g_i) = \pi_2(h_i), i = 1, 2$. Thus $(g_i, h_i) \in \widetilde{\Lieg}, i = 1, 2$, and we obtain
\[
\begin{aligned}
& [(g_1,h_1),(g_2,h_2)]+[(g_2,h_2),(g_1,h_1)] =\\
 &([g_1,g_2]+[g_2,g_1], [h_1,h_2]+[h_2,h_1]) =\\
&([g_1,g_2]+[g_2,g_1], \xi([g_1,g_2]+[g_2,g_1])).
\end{aligned} \]
This implies that $[\widetilde{\Lieg}, \widetilde{\Lieg}]_{\Lie} = \{ (g', \xi(g') \mid g' \in [\Lieg_1, \Lieg_1]_{\Lie}\}$.

  For the case of $i=1,$ $\gamma_1=Id_{\Lieq_1}.$ Now let $(g',\xi(g'))\in \Ker(\tau_1) \cap[\widetilde{\Lieg},\widetilde{\Lieg}]_{\Lie}.$ Then $\tau_1(g', \xi(g'))= g' =0$, hence $\xi(g')=0$  and thus $\Ker(\tau_1) \cap[\widetilde{\Lieg},\widetilde{\Lieg}]_{\Lie} = 0$. Proposition \ref{equivalence} implies that $(g_1) \sim (\widetilde{g})$.

   For the case of $i=2,$ $\gamma_2=\eta$ and the fact $\Ker(\tau_2)\cap[\widetilde{\Lieg},\widetilde{\Lieg}]_{\Lie}=0$ can be proved in an analogous way as case $i=1$. Proposition \ref{equivalence} implies that $(g_2) \sim (\widetilde{g})$.  \rdg
\bigskip

Let $(g) :  0 \to \Lien \stackrel{\chi}\to  \Lieg \stackrel{\pi} \to \Lieq \to 0$ be a $\Lie$-central extension and  $\Liea:=\Lieh/[\Lieh,\Lieh]_{\Lie}$ where $\Lieh$ is a Leibniz algebra. Let $\varphi : \Lieg \times \Liea \to \Lieg$ be  the projection onto $\Lieg$, then
\begin{equation}
(g \times a) :  0 \to \Lien \times \Liea  \stackrel{\chi \times id}\to  \Lieg \times \Liea \stackrel{\pi \circ \varphi} \to \Lieq \to 0
\end{equation}
is a $\Lie$-central extension.

Let $\mu : \Lieg \to \Lieg \times \Liea, \mu(g)=(g,0)$; $\varphi' : \Lien \times \Liea \to \Lien, \varphi'(n,a)=n$; $\mu' : \Lien \to \Lien \times \Liea, \mu'(n)=(n,0)$ be. Then the diagrams:
\[
\xymatrix{
(g \times a) :  0 \ar[r] & \Lien  \times \Liea \ar[r]^{\chi \times id} \ar@<0.5ex>[d]^{\varphi'} & \Lieg \times \Liea \ar[r]^{\pi \circ \varphi} \ar@<0.5ex>[d]^{\varphi} &  \Lieq \ar@{=}[d] \ar[r]& 0 \\
(g) :  0 \ar[r] & \Lien \ar[r]^{\chi} \ar@<0.5ex>[u]^{\mu'} & \Lieg \ar[r]^{\pi} \ar@<0.5ex>[u]^{\mu} & \Lieq \ar[r] & 0
}
\]
are commutative, $(\varphi', \varphi, id_{\Lieq}) : (g \times a) \to (g)$ is a $\Lie$-isoclinic epimorphism and $(\mu', \mu, id_{\Lieq}) : (g) \to (g \times a)$ is a $\Lie$-isoclinic monomorphism.

Let $(g_i), \sigma_i, \tau_i, \gamma_i, i=1,2,$ and $(\widetilde{g})$ be as in (\ref{tilda g}), $\eta : \Lieq_1 \to \Lieq_2$ an isomorphism and $\Liea = \widetilde{\Lieg}/[\widetilde{\Lieg}, \widetilde{\Lieg}]_{\Lie}$, and $\alpha : \Lien_1  \times \Lien_2 \to \Lien_1  \times \Liea, \alpha(n) = (\sigma_1(n), \overline{\lambda(n)})$; $\beta : \widetilde{\Lieg} \to \Lieg_1 \times \Liea, \beta(g) = (\tau_1(g), \overline{g})$.

\begin{Le} \label{monomorphism}
Assume that $\eta$ induces a $\Lie$-isoclinism from $(g_1)$ to $(g_2)$. Then $(\alpha, \beta, \gamma_1) : (\widetilde{g}) \to (g_1 \times a)$ is a $\Lie$-isoclinic monomorphism.
\end{Le}
{\it Proof.} Let  $g\in \Ker(\beta) \cap[\widetilde{\Lieg},\widetilde{\Lieg}]_{\Lie}.$ Then $\beta(g)=0$ i.e. $(\tau_1(g), \bar{g}) =0$  which implies that $\tau_1(g)=0$ and  $\bar{g}=0.$ So  $g\in \Ker(\tau_1) \cap[\widetilde{\Lieg},\widetilde{\Lieg}]_{\Lie} = 0$ by Proposition 3.14. Hence $g=0$ and thus $\Ker(\beta) \cap[\widetilde{\Lieg},\widetilde{\Lieg}]_{\Lie}=0.$ To show that $\beta$ is one-to-one, notice that  $\bar{g}=0$ iff $g\in[\widetilde{\Lieg},\widetilde{\Lieg}]_{\Lie}.$ So $\beta(g)=0$ iff  $g\in \Ker(\tau_1) \cap[\widetilde{\Lieg},\widetilde{\Lieg}]_{\Lie} = 0,$ and thus $g=0.$  \rdg

\bigskip

Let $(\eta, \xi) : (g_1) \sim (g_2)$ be a $\Lie$-isoclinism. Now we consider the following commutative diagram where $\Liea:=\widetilde{\Lieg}/[\widetilde{\Lieg},\widetilde{\Lieg}]_{\Lie}$:

\begin{equation} \label{diagram}
\xymatrix{
(g_2) : & 0 \ar[r] & \Lien_2 \ar[r]^{\chi_2} & \Lieg_2 \ar[r]^{\pi_2} & \Lieq_2 \ar[r] & 0\\
(\widetilde{g}) : & 0 \ar[r] & \Lien_1 \times \Lien_2 \ar[r]^{\lambda} \ar@{>>}[u]_{\sigma_2} \ar@{>->}[d]^{\alpha} \ar@/_6mm/@{-->}[dd]& \widetilde{\Lieg} \ar[r]^{\rho} \ar@{>>}[u]_{\tau_2} \ar@{>->}[d]^{\beta} \ar@/_6mm/@{-->}[dd]  & \Lieq \ar[r] \ar[u]_{\gamma_2}^{\wr} \ar[d]^{\gamma_1}_{\wr}& 0\\
(g_1 \times \frak{a}) : & 0 \ar[r] & \Lien_1 \times \frak{a} \ar[r] \ar@{>>}[d]^{nat'}& \Lieg_1 \times \frak{a} \ar[r] \ar@{>>}[d]^{nat}& \Lieq_1 \ar[r] \ar@{=}[d]& 0\\
nat' (g_1 \times \frak{a}) = \left(\frac{g_1 \times \frak{a}}{\alpha(\Lien_1)}\right) : & 0 \ar[r] & \frac{\Lien_1 \times \frak{a}}{\alpha(\Lien_1)}\quad \ar@{>->}[r] & \frac{\Lieg_1 \times \frak{a}}{\beta \lambda(\Lien_1)} \ar@{>>}[r] & \Lieq_1 \ar[r] & 0\\
(g_1 \times \frak{a}) : & 0 \ar[r] & \Lien_1 \times \frak{a} \quad \ar@{>->}[r] \ar@{>>}[u]_ {nat'}& \Lieg_1 \times \frak{a} \ar@{>>}[r] \ar@{>>}[u]_{nat}& \Lieq_1 \ar[r] \ar@{=}[u]& 0\\
(g_1) : & 0 \ar[r] & \Lien_1 \quad \ar@{>->}[r]^{\chi_1} \ar[u]_{\mu'} \ar@/^8mm/@{>-->}[uu]& \Lieg_1 \ar@{>>}[r]^{\pi_1} \ar[u]_{\mu} \ar@/^6mm/@{>-->}[uu] & \Lieq_1 \ar[r] \ar@{=}[u]& 0
}
\end{equation}

The following holds:
\begin{equation}
\alpha(\Lien_1) = \{ (n_1, (\overline{\chi_1(n_1),0)}\ ), {\rm for \ all \ }  n_1 \in \Lien_1\}
\end{equation}
\begin{equation}
\Ker(nat) = \{ (\chi_1(n_1), \overline{(\chi_1(n_1), 0)}\ ), {\rm for \ all \ }  n_1 \in \Lien_1\}
\end{equation}
\begin{equation}
[\widetilde{\Lieg}, \widetilde{\Lieg}]_{\Lie} = \{[g, \xi(g)],  {\rm for \ all \ }  g \in [\Lieg_1, \Lieg_1]_{\Lie}\}.
\end{equation}
This yields:
\begin{equation} \label{intersection}
\Ker(nat) \cap [\Lieg_1 \times \frak{a}, \Lieg_1 \times \frak{a}]_{\Lie} = 0
\end{equation}
\begin{equation}
\alpha(\Lien_1) \cap \mu'(\Lien_1) =0.
\end{equation}
Furthermore we have
\begin{equation}
\Ker(nat' \circ \alpha) = \Lien_1 = \Ker(\sigma_2)
\end{equation}
and we obtain:

\begin{Le}\label{lemma natural}\
\begin{enumerate}
\item[a)] The composition of $(nat', nat, Id_{\Lieq_1})$ and $(\mu', \mu, Id_{\Lieq_1})$ is a \Lie-isoclinic monomorphism from $(g_1)$ into $nat'((g_1 \times \frak{a}))$.

    \item[b)] The composition of $(nat', nat, Id_{\Lieq_1})$ and $(\alpha, \beta, \gamma_1)$ induces a \Lie-isoclinic monomorphism from $(g_2)$ into $nat' ((g_1 \times \frak{a}))$.
\end{enumerate}
\end{Le}
{\it Proof.}
{\it a)} Let $g\in \Ker(nat\circ \mu)\cap [\Lieg_1,\Lieg_1]_{\Lie}.$ Then $\mu(g)\in \Ker(nat)$ and $g=[g_1,g_2]+[g_2,g_1]$ for some $g_1,g_2\in \Lieg_1.$ So $\mu(g)=[\mu(g_1),\mu(g_2)]+[\mu(g_2),\mu(g_1)]\in [\Lieg_1 \times \frak{a}, \Lieg_1 \times \frak{a}]_{\Lie}.$ Thus $(g,0)=\mu(g)\in \Ker(nat) \cap [\Lieg_1 \times \frak{a}, \Lieg_1 \times \frak{a}]_{\Lie} = 0$ by (\ref{intersection}). Hence  the composition  $(nat', nat, Id_{\Lieq_1})$ and $(\mu', \mu, Id_{\Lieq_1})$ is a \Lie-isoclinism by Proposition \ref{equivalence}.

To show that $nat\circ\mu$ is one-to-one, let $g\in\Lieg_1$ such that $(nat\circ\mu)(g)=0.$ This implies that $(g,0)\in \beta \lambda(\Lien_1)=(\chi_1(n),\overline{(\chi_1(n),0)})$ for some $n\in \Lien_1.$ So $g = \chi_1(n)$ and $(\chi_1(n),0)\in [\widetilde{\Lieg},\widetilde{\Lieg}]_{\Lie}.$ Therefore $(g,0)\in \Ker(\tau_2)\cap [\widetilde{\Lieg},\widetilde{\Lieg}]_{\Lie}.$ Since  $(\eta,\xi)$ is  a $\Lie$-isoclinism from $(g_1)$ to $(g_2),$ it follows by the proof of Proposition \ref{Lie iso epi} that $\Ker(\tau_2)\cap [\widetilde{\Lieg},\widetilde{\Lieg}]_{\Lie}=0.$  Hence $g=0.$
\bigskip

{\it b)}  Consider the map $\delta:\Lieg_2\to \frac{\Lieg_1 \times \frak{a}}{\beta \lambda(\Lien_1)}$ defined by $\delta(h)=(nat\circ\beta)(x,h)$ where $\pi_2(h)=\eta(\pi_1(x))$ for some $x\in\Lieg_1.$

$\delta$ is well defined, since for $x,x'\in\Lieg_1$ such that $\eta(\pi_1(x))=\pi_2(h)=\eta(\pi_1(x')), $ then $x-x'\in \Ker(\pi_1)$, that is $x-x'=\chi_1(n)$ for some $n\in\Lieg_1.$ We then have
 $$(x,\overline{(x,h)})-(x',\overline{(x',h)})=(x-x',\overline{(x-x',0)})=(\chi_1(n),\overline{(\chi_1(n),0)})\in \beta \lambda(\Lien_1).$$  Hence $(nat\circ\beta)(x,h)=(nat\circ\beta)(x',h).$

 Let $\delta_0$ be the restriction of $\delta$ to $\Lien_2.$  It is easy to check that $(\delta_0,\delta,\eta^{-1}): (g_2)\to nat' (g_1 \times \frak{a})$ is a homomorphism of \Lie-central extensions.

 Let $h\in \Ker(\delta).$ Then $(x,\overline{(x,h)})\in \beta \lambda(\Lien_1)$ for some $x\in\Lieg_1$ satisfying $\pi_2(h)=\eta(\pi_1(x)).$  Then $(x,\overline{(x,h)})=(\chi_1(n),\overline{(\chi_1(n),0)})$ for some $n\in \Lien_1.$ So $x=\chi_1(n)$ and $(0,h)=(x-\chi_1(n),h)\in[\widetilde\Lieg,\widetilde\Lieg]_{\Lie}.$ So $(0,h)\in \Ker(\tau_1)\cap[\widetilde\Lieg,\widetilde\Lieg]_{\Lie}=0$ by the proof of Proposition 3.14. Therefore $h=0.$ So $\delta$ is one to one. This implies that  $\Ker(\delta)\cap [\Lieg_2,\Lieg_2]_{\Lie}=0.$ It  follows
 by Proposition 3.10 that $(\delta_0,\delta,\eta^{-1})$ is a $\Lie$- isoclinic monomorphism. \rdg

\begin{Th} \label{isoclinic epi}
The following statements are equivalent:
\begin{enumerate}
\item[a)] The \Lie-central extensions $(g_1)$ and $(g_2)$ are \Lie-isoclinic.

\item[b)] There exists a \Lie-central extension $(g')$ together with \Lie-isoclinic epimorphism from $(g')$ onto $(g_1)$ and $(g_2)$.

\item[c)] There exists a \Lie-central extension $(g'')$ together with \Lie-isoclinic monomorphisms from $(g_1)$ and $(g_2)$ into  $(g'')$.
\end{enumerate}
\end{Th}
{\it Proof.} {\it a)} $\iff$ {\it b)} Apply Proposition \ref{Lie iso epi} {\it b)} by letting $(g')=(\tilde{g})$.

\noindent {\it a)} $\implies$ {\it c)} Apply Lemma \ref{lemma natural} by letting $(g'')=(nat'(g_1\times\Liea)).$

\noindent  {\it c)} $\implies$ {\it a)} Direct consequence of Proposition \ref{equivalence relation}.
 \rdg

\begin{Rem}
If $(g_1)$ and $(g_2)$ are finite-dimensional, extensions $(g')$ and $(g'')$ in Theorem \ref{isoclinic epi} can also be chosen finite-dimensional. The situation in Theorem \ref{isoclinic epi} can be roughly sketched by $(g_1) \twoheadleftarrow (g')  \twoheadrightarrow (g_2)$ and $(g_1) \rightarrowtail (g'') \leftarrowtail (g_2)$. In particular, the equivalence relations for \Lie-central extensions generated by \Lie-isoclinic epimorphisms, respectively \Lie-isoclinic monomorphisms, coincide with \Lie-isoclinism.
\end{Rem}

\begin{Th} \label{iso mono}
The following statements are equivalent:
\begin{enumerate}
\item[a)] $(g_1)$ and $(g_2)$ are \Lie-isoclinic.

\item[b)] There exist a Leibniz algebra with trivial \Lie-commutator $\frak{a}$, a \Lie-central extension $(g')$, a \Lie-isoclinic monomorphism from  $(g')$ into $(g_1 \times \frak{a})$ and a \Lie-isoclinic epimorphism from $(g')$ onto $(g_2)$.

    \item[c)] There exist a Leibniz algebra with trivial \Lie-commutator $\frak{b}$, a \Lie-central extension $(g'')$, a \Lie-isoclinic epimorphism from $(g_1 \times \frak{b})$ onto $(g'')$, and a \Lie-isoclinic monomorphism from $(g_2)$ into $(g'')$.
\end{enumerate}
\end{Th}
{\it Proof.} {\it a)} $\implies$ {\it b)} It follows from Proposition \ref{Lie iso epi} and Lemmas \ref{monomorphism} and \ref{lemma natural} by letting $\Liea=\tilde{\Lieg}/[\tilde{\Lieg},\tilde{\Lieg}]_{\Lie}$ and $(g')=(\tilde{\Lieg}).$

\noindent {\it b)} $\implies$ {\it c)} By letting $\Lieb=\tilde{\Lieg}/[\tilde{\Lieg},\tilde{\Lieg}]_{\Lie}$ and $(g'')=(nat'(g_1\times\Liea)),$ the implication  follows by definition of $(nat'(g_1\times\Liea))$ and the proof of Lemma \ref{lemma natural} {\it b)}.

\noindent {\it c)} $\implies$ {\it a)} It suffices to have in mind the following chain of \Lie-isoclinisms:
\[ \xymatrix{
(g_1)\quad \ar@{>->} @<3pt>[r] &(g_1 \times \Lieb) \ar@{>>}[r]  \ar@{>>} @<3pt>[l] & (g'')& \quad(g_2) \ar@{>->}[l]}.
\] \rdg

\begin{Pro} \label{natural}
Let $\Lieg$ be a Leibniz algebra.
\begin{enumerate}
\item[a)]  Let $\frak{a}$ be a Leibniz algebra with trivial \Lie-commutator. Then $\Lieg$ and $\Lieg \times \frak{a}$ are \Lie-isoclinic.

\item[b)] Let $\Lien$ be a two-sided ideal of \Lieg. The natural homomorphism $nat : \Lieg \twoheadrightarrow \Lieg/\Lien$ is a \Lie-isoclinic epimorphism if and only if $\Lien \cap [\Lieg, \Lieg]_{\Lie} = 0$.

 \item[c)] Let $\Lieh$ be a subalgebra of $\Lieg$. The embedding of $\Lieh$ into $\Lieg$ is a \Lie-isoclinic monomorphism if and only $\Lieh + Z_{\Lie}(\Lieg) = \Lieg$.
\end{enumerate}
\end{Pro}
{\it Proof.} {\it a)} Notice that for a Leibniz algebra $\Lieg$ and a Leibniz algebra with trivial \Lie-commutator $\Liea,$ we have $Z_{\Lie}(\Lieg\times\Liea) = Z_{\Lie}(\Lieg)\times\Liea.$

The canonical inclusion $\beta : \Lieg \to \Lieg \times \Liea, \beta(g)=(g,0)$, satisfies that $\Ker(\beta) \cap [\Lieg, \Lieg]_{\Lie} = 0$ and any element $(g, a) \in \Lieg \times \Liea$ can be written as  $(g,a) = (g,0)+ (0,a) \in \Image(\beta)+ Z_{\Lie}(\Lieg \times \Liea)$. Now Proposition \ref{character}  ends the proof.

{\it b)} Follows by Proposition \ref{character} since $\Ker(nat)=\Lien.$ In addition, if  $\Lien \cap [\Lieg, \Lieg]_{\Lie} \neq 0,$ then $nat$ is not a \Lie-isoclinic homomorphism by Proposition \ref{equivalence}.

Conversely, if $\Lien \cap [\Lieg, \Lieg]_{\Lie}=0,$ then $\xi:[\Lieg,\Lieg]_{\Lie}\to[\Lieg/\Lien,\Lieg/\Lien]_{\Lie}, \xi([g_1,g_2]+[g_2,g_1])=[g_1+\Lien,g_2+\Lien]+[g_2+\Lien,g_1+\Lien]$, is an  isomorphism and $Z_{\Lie}(\Lieg/\Lien)=\frac{Z_{\Lie}(\Lieg)}{\Lien},$ and $\frac{\Lieg}{Z_{\Lie}(\Lieg)}\stackrel{\eta}{\cong} \frac{\frac{\Lieg}{\Lien}}{\frac{Z_{\Lie}(\Lieg)}{\Lien}}$ by the third isomorphism theorem. Now the commutativity of diagram (\ref{square isoclinic})  immediately follows.

{\it c)} Since $\Image(\Lieh\hookrightarrow\Lieg)=\Lieh$ and $\Ker(\Lieh\hookrightarrow\Lieg)=0$, then  Proposition \ref{character} ends the proof.
 \rdg

\begin{Th}
Let $\Lieg$ and $\Lieq$ be Leibniz algebras. Then the following properties are equivalent:
\begin{enumerate}
\item[a)] $\Lieg$ and $\Lieq$ are \Lie-isoclinic.

\item[b)] There exist a Leibniz algebra with trivial \Lie-commutator $\frak{a}$, a subalgebra $\Lieh$ of $\Lieg \times \frak{a}$ with $\Lieh + Z_{\Lie}(\Lieg \times \frak{a}) = \Lieg \times \frak{a}$, and a two-sided ideal $\Lien$ of $\Lieh$ with $\Lien \cap [\Lieh, \Lieh]_{\Lie} =0$, such that $\Lieh/\Lien$ is isomorphic to $\Lieq$.

    \item[c)] There exist a Leibniz algebra with trivial \Lie-commutator $\Lieb$, a two-sided ideal $\Liem$ of $\Lieg \times \Lieb$ with $\Liem \cap [\Lieg \times \Lieb, \Lieg \times \Lieb]_{\Lie} = 0$, and a subalgebra $\Liei$ of $\frac{\Lieg \times \Lieb}{\Liem}$, with $\Liei + Z_{\Lie}(\frac{\Lieg \times \Lieb}{\Liem}) = \frac{\Lieg \times \Lieb}{\Liem}$, such that $\Liei$ is isomorphic to $\Lieq$.
\end{enumerate}
\end{Th}
{\it Proof.} {\it a)} $\iff$ {\it b)} Assume a), then by Theorem \ref{iso mono}  there exists a Leibniz algebra with trivial \Lie-commutator $\frak{a}$, a \Lie-central extension $(g')$, a \Lie-isoclinic monomorphism (denoted by $\alpha$) from  $(g')$ into $(g \times \frak{a})$ and a \Lie-isoclinic epimorphism  (denoted by $\beta$) from $(g')$ onto $(q)$.  Let $\Lieh=\Image(\alpha)\cong \Lieg'$ and $\Lien:=\Ker(\beta).$ Then $\Lieh + Z_{\Lie}(\Lieg \times \Liea) = \Lieg \times \Liea$  by Proposition \ref{natural}, $\frac{\Lieg'}{\Lien}\cong\Lieq$ and we have by Proposition \ref{character}  that $\Lien \cap [\Lieh, \Lieh] = 0.$

Conversely, by letting $\Lieg'=\Lieh$ then Proposition \ref{natural} provides a \Lie-isoclinic monomorphism  from  $(g')$ into $(g \times \frak{a})$ and a \Lie-isoclinic epimorphism   from $(g')$ onto $(q)$. The result follows by Theorem \ref{iso mono}

{\it a)} $\iff${\it c)} Assume a), then by Theorem \ref{iso mono} there exists a  Leibniz algebra with trivial \Lie-commutator $\frak{b}$, a central extension $(g'')$, a \Lie-isoclinic epimorphism (denoted by $\alpha$) from $(g\times \frak{b})$ onto $(g'')$, and a \Lie-isoclinic monomorphism (denoted by $\beta$) from $(q)$ into $(g'')$. Let $\Liem=\Ker(\alpha).$ Then $\frac{g\times \frak{b}}{\Liem}\cong \Lieg''$ and  $\Liem \cap [\Lieg \times \Lieb, \Lieg \times \Lieb]_{\Lie} = 0$ by Proposition \ref{character}. Also by Proposition \ref{natural}  we have $\Liei+Z_{\Lie}(\Lieg'')=\Lieg'',$ where $\Liei = \Image(\beta) \cong \Lieq$.

Conversely, by letting $\Lieg''=\frac{g\times \frak{b}}{\Liem}$ we have by Proposition \ref{natural} that the embedding $\Lieq\cong \Liei$ into $\Lieg''$ is a \Lie-isoclinic monomorphism. Also  Proposition \ref{natural} implies a \Lie-isoclinic epimorphism from $\Lieg\times \Lieb$ to $\Lieg''.$

The result follows by Theorem \ref{iso mono}. \rdg

\begin{Rem}
Let $\beta : \Lieg \to \Lieq$ be a \Lie-isoclinic epimorphism. Then $\beta$ induces an isoclinism $(\eta, \xi) : \Lieg \sim \Lieq$, where $\xi = \beta_{\mid [\Lieg, \Lieg]_{\Lie}}$. Hence we have the following commutative diagram:
\begin{equation} \label{commutative diagram}
\xymatrix{
0 \ar[r] & [\Lieq, \Lieq]_{\Lie} \ar[r]^{\quad \xi^{-1}} \ar@{=}[d] & \Lieg \ar[r] \ar@{>>}[d]^{\beta} & \frak{a} \ar[r] \ar@{>>}[d] & 0\\
0 \ar[r] & [\Lieq, \Lieq]_{\Lie} \ar@{^{(}->}[r]  & \Lieq \ar[r]  & \frac{\Lieq}{[\Lieq, \Lieq]_{\Lie}} \ar[r]  & 0
} \end{equation}
where $\frak{a} \cong \frac{\Lieg}{[\Lieg, \Lieg]_{\Lie}}$. On the other hand, if we fix $\Lieq$, while $\frak{a}$ is a Leibniz algebra with trivial \Lie-commutator and $\Lieg$ a Leibniz algebra which fits into (\ref{commutative diagram}) for some epimorphism $\beta$, then $[\Lieg, \Lieg]_{\Lie} = \Image(\xi^{-1})$ and $\Ker(\beta) \cap [\Lieg, \Lieg]_{\Lie} =0$. Hence, diagram (\ref{commutative diagram}) determines all groups which map epi-isoclinic onto $\Lieq$.
\end{Rem}


\section{\Lie-Isoclinism and the Schur \Lie-multiplier} \label{Schur}

In this section we analyze the connection between \Lie-isoclinism and the second \Lie-homology with trivial coefficients, which we call  the Schur \Lie-multiplier thanks to the isomorphism (\ref{Hopf}).

Consider the $\Lie$-central extensions $(g) : 0 \to \Lien \stackrel{\chi} \to \Lieg \stackrel{\pi} \to \Lieq \to 0$ and $(g_i) : 0 \to \Lien_i \stackrel{\chi_i}\to \Lieg_i \stackrel{\pi_i} \to \Lieq_i \to 0, i=1, 2$. By Proposition 2 in \cite{CK},  associated to $(g)$ there exists the following  six-term exact sequence
\begin{equation}\label{six-term}
 \Lien\otimes \Lieg_{\Lie} \longrightarrow \HL^{\Lie}_2(\Lieg) \longrightarrow \HL^{\Lie}_2(\q) \stackrel{\theta(g)}\longrightarrow \frak{n} \longrightarrow \HL^{\Lie}_1(\Lieg) \longrightarrow \HL^{\Lie}_1(\q) \longrightarrow 0.
 \end{equation}
 and from its proof easily follows that
 \begin{equation} \label{image}
 \Image(\theta(g)) = \Lien \cap [\Lieg, \Lieg]_{\Lie}.
 \end{equation}
 Thus $\theta(g)$ induces a homomorphism $\theta'(g) : \HL_2^{\Lie}(\Lieq) \to [\Lieg, \Lieg]_{\Lie}$ such that the sequence
 \begin{equation} \label{four term}
0 \to  \Ker(\theta(g)) \to \HL_2^{\Lie}(\Lieq) \stackrel{\theta'(g)}\to [\Lieg, \Lieg]_{\Lie} \stackrel{\pi'} \to [\Lieq, \Lieq]_{\Lie} \to 0
 \end{equation}
 is exact, where $\pi'$ is induced by $\pi : \Lieg \twoheadrightarrow \Lieq$. The naturality of  $\theta(g)$ with respect to homomorphisms of \Lie-central extensions (see \cite[Theorem 5.9]{EVDL}) yields the naturality of the sequence (\ref{four term}).

For the $\Lie$-central extension $(g)$, $\theta(g)$ is a monomorphism if and only if $\Image(\HL_2^{\Lie}(\pi))$ $= 0$  by the exactness of (\ref{six-term}), equivalently $\HL_2^{\Lie}(\pi)$ is the zero map, or equivalently, the sequence $0 \to \HL^{\Lie}_2(\q) \stackrel{\theta(g)}\to \frak{n} \to \HL^{\Lie}_1(\Lieg) \to \HL^{\Lie}_1(\q) \to 0$ is exact.

On the other hand, by (\ref{image})  $\theta(g)$ is a monomorphism if and only if $\HL^{\Lie}_2(\q) \cong \Lien \cap [\Lieg, \Lieg]_{\Lie}$. Having in mind \cite[Proposition 3]{CK}, the \Lie-central extension $(g)$ with $\theta(g)$  a monomorphism is a \Lie-trivial extension whenever $\HL^{\Lie}_2(\q) = 0$, so $\HL^{\Lie}_2(\q)$ measures the deficiency of  a \Lie-central extension to be a \Lie-trivial extension. For this reason, we call a such \Lie-central extension a {\it quasi \Lie-trivial extension}.

 \begin{Le} \label{A}
Let $(\alpha, \beta, \gamma):(g_1) \to (g_2)$ be a homomorphism of \Lie-central extensions and assume that $\gamma : \Lieq_1 \to \Lieq_2$ is an isomorphism. Then the following statements are equivalent:
\begin{enumerate}
\item[a)] $(\alpha, \beta, \gamma)$ is a \Lie-isoclinism.
\item[b)] $\HL_2^{\Lie}(\gamma) \left( \Ker(\theta(g_1)) \right) = \Ker(\theta(g_2))$.
\end{enumerate}
 \end{Le}
{\it Proof.} The naturality of (\ref{four term}) induces the following commutative diagram:
\[
\xymatrix{
 & & \Lien_1 \cap [\Lieg_1,\Lieg_1]_{\Lie} \ar@{>->}[dr]& & & \\
 \Ker(\theta(g_1)) \ar[d] \ar@{>->}[r] & \HL_2^{\Lie}(\Lieq_1) \ar[rr]^{\theta'(g_1)} \ar[d]^{\HL_2^{\Lie}(\gamma)} \ar@{>>}[ur]^{\theta(g_1)} &  & [\Lieg_1,\Lieg_1]_{\Lie} \ar@{>>}[r]^{\pi_1'} \ar[d]^{\beta'} & [\Lieq_1,\Lieq_1]_{\Lie} \ar[d]^{\gamma'}\\
  \Ker(\theta(g_2)) \ar@{>->}[r] & \HL_2^{\Lie}(\Lieq_2) \ar[rr]^{\theta'(g_2)}  \ar@{>>}[dr]^{\theta(g_2)} &  & [\Lieg_2,\Lieg_2]_{\Lie} \ar@{>>}[r]^{\pi_2'}  & [\Lieq_2,\Lieq_2]_{\Lie}\\
 & & \Lien_2 \cap [\Lieg_2,\Lieg_2]_{\Lie} \ar@{>->}[ur]& &  &
}
\]
where $\theta'(g_i), i=1,2$, are defined as above, $\beta' = \beta_{\mid [\Lieg_1,\Lieg_1]_{\Lie}}, \gamma' = \gamma_{\mid [\Lieq_1,\Lieq_1]_{\Lie}}$.

Since $\gamma$ is a isomorphism, then $\gamma'$ and $\HL_2^{\Lie}(\gamma)$ are also isomorphisms, the restriction of $\HL_2^{\Lie}(\gamma)$ to $\Ker(\theta(g_1))$ is a monomorphism and $\beta$ is an epimorphism.

By the commutativity of the left hand square in the diagram, we have that $\HL_2^{\Lie}(\gamma)(\Ker(\theta(g_1))) \subseteq \Ker(\theta (g_2))$. Conversely, for any $y \in \Ker(\theta (g_2))$, there exists $x \in  \HL_2^{\Lie}(\Lieq_1)$ such that $\HL_2^{\Lie}(\gamma)(x) =y$. Now $0= \theta'(g_2)(y) = \theta'(g_2) \circ \HL_2^{\Lie}(\gamma)(x) = \beta'_{\mid} \circ \theta'(g_1)(x)$. Hence $x \in \Ker(\theta(g_1))$ whenever $\beta'_{\mid}$ is a monomorphism.

Consequently, statement {\it b)} holds if and only if $\beta'_{\mid}$ is a monomorphism if and only if $(\alpha, \beta, \gamma)$ is a \Lie-isoclinism (having in mind Proposition \ref{equivalence}). \rdg

\begin{Le} \label{B}
Let $\eta : \Lieq_1 \to \Lieq_2$ be an isomorphism and $(\widetilde{g}) = (g_1 \times g_2^{\eta})^{\Delta}$ as in section \ref{Lie iso Lb alg}, then $$\Ker(\theta(\widetilde{g})) =\HL_2^{\Lie}(\eta)^{-1}(\Ker(\theta(g_2))) \cap \Ker(\theta(g_1)).$$
\end{Le}
{\it Proof.} From diagram (\ref{tilda g})  and the naturality of $\theta(g)$  we have that
$$\theta(\widetilde{g}) = \theta (g_1) \times \theta(g_2) \circ \HL_2^{\Lie}(\gamma_2).$$
Since $\HL_2^{\Lie}(\gamma_2)$ is an isomorphism, then $\Ker(\theta(g_2) \circ \HL_2^{\Lie}(\gamma_2)) = \HL_2^{\Lie}(\gamma_2)^{-1}(\Ker(\theta(g_2)))$, hence the required equality follows. \rdg

\begin{Th} \label{main1}
Let $\eta : \Lieq_1 \to \Lieq_2$ be an isomorphism. Then the following statements are equivalent:
\begin{enumerate}
\item[a)] $\eta$ induces a \Lie-isoclinism from $(g_1)$ to $(g_2)$.
\item[b)] There exists an isomorphism $\beta' : [\Lieg_1, \Lieg_1]_{\Lie} \to [\Lieg_2, \Lieg_2]_{\Lie}$ with $\beta' \circ \theta'(g_1) = \theta'(g_2) \circ \HL_2^{\Lie}(\eta)$.
  \item[c)] $\HL_2^{\Lie}(\eta)(\Ker(\theta(g_1))) = \Ker(\theta(g_2))$.
\end{enumerate}
\end{Th}
{\it Proof.} {\it a)} $\Rightarrow$ {\it b)} Let $(\eta, \beta') : (g_1)\sim (g_2)$ be. Then Propositions \ref{Lie iso epi} and \ref{equivalence} imply that $\beta' = \tau_2' \circ \tau_1'^{-1}$, where $\tau_i' = \tau_i{_{\mid [\widetilde{\Lieg}, \widetilde{\Lieg}]_{\Lie}}}$ and the proof of Proposition \ref{Lie iso epi} shows that both of them are isomorphisms.

Then the naturality of sequence (\ref{four term}) applied to diagrams (\ref{tilda g}) implies
$$\theta'(g_2) \circ \HL_2^{\Lie}(\gamma_2) = \tau'_2 \circ \theta'(\widetilde{\Lieg})= \tau'_2 \circ \tau_1'^{-1} \circ \theta'(g_1) \circ \HL_2^{\Lie}(\gamma_1),$$
hence the required equality.

{\it b)} $\Rightarrow$ {\it c)} Since $\beta' \circ \theta'(g_1) = \theta'(g_2) \circ \HL_2^{\Lie}(\eta)$, then $\HL_2^{\Lie}(\eta) (\Ker(\theta(g_1))) \subseteq \Ker(\theta(g_2))$. The converse inclusion is followed  thanks to be $\beta'$ an isomorphism.

{\it c)} $\Rightarrow$ {\it a)} Let $(\widetilde{g})$ be as above, then Lemma \ref{B} implies that
\begin{equation} \label{equality}
\Ker(\theta(\widetilde{g})) = \HL_2^{\Lie}(\eta)^{-1} \circ \HL_2^{\Lie}(\eta) (\Ker(\theta(g_1))) \cap \Ker(\theta(g_1)) = \Ker(\theta(g_1)).
\end{equation}
From diagram (\ref{tilda g}), the epimorphisms $(\sigma_i, \tau_i, \gamma_i) : (\widetilde{g}) \sim (g_i)$ and the equality (\ref{equality}) imply that
$$\HL_2^{\Lie}(\gamma_i)(\Ker(\theta(\widetilde{g}))) = \Ker(\theta(g_i)), i = 1, 2.$$
Lemma \ref{A} implies that $(\sigma_i, \tau_i, \gamma_i) : (\widetilde{g}) \sim (g_i), i = 1, 2.$ Now Proposition \ref{Lie iso epi} ends the proof. \rdg

\begin{Co}\
\begin{enumerate}
\item[a)] Let $Autcl(g)$ be the group of \Lie-autoclinisms of $(g)$. Then
$$Autcl(g) \cong \{ \eta \mid \eta \in Aut(\Lieq), {\rm such\ that}\ \HL_2^{\Lie}(\Ker(\theta(g))) = \Ker(\theta(g))\}.$$
\item[b)] If $(g)$ is a quasi \Lie-trivial extension, then $Autcl(g) \cong Aut(\Lieq)$.
\item[c)] Any two quasi \Lie-trivial extensions of a given Leibniz algebra are \Lie-isoclinic.
\end{enumerate}
\end{Co}
{\it Proof.} Direct application of Theorem \ref{main1}. \rdg
\bigskip

Following \cite{CK}, a  $\Lie$-central extension $0 \to \Lien \stackrel{\chi} \to  \Lieg \stackrel{\pi}\to \Lieq \to 0$ is called a $\Lie$-stem extension if $\Lieg_{\Lie} \cong \Lieq_{\Lie}$. A $\Lie$-stem cover is a $\Lie$-stem extension such that the induced map $HL^{\Lie}_2(\Lieg) \to HL^{\Lie}_2(\Lieq)$ is the zero map. Proposition 6 in \cite{CK} characterizes $\Lie$-stem covers  by the fact that the induced map $\theta : HL_2^{\Lie}(\Lieq) \to \Lien$ is an isomorphism.

Consider $\Lieq$ a Leibniz algebra for which there exists a free presentation $0 \to \Lies \to \Lief \stackrel{\tau}\to \Lieq \to 0$ satisfying $\Lies \subseteq [\Lief, \Lief]_{\Lie}$. For instance $\Lieq$ is a free Leibniz algebra, then $\Lief \cong \Lieq, \Lies=0$.
Under this assumption, it is easy to check that the map $HL_2^{\Lie}(\Lieq) \cong \frac{\Lies \cap [\Lief, \Lief]_{\Lie}}{[\Lies, \Lief]_{\Lie}} \to \frac{\Lies}{[\Lies, \Lief]_{\Lie}}$ is an epimorphism, so Proposition 5 (c) and Example 2 (b) in \cite{CK} provides the following $\Lie$-stem cover:
\begin{equation} \label{Lie stem cover}
0 \to \frac{\Lies}{[\Lies, \Lief]_{\Lie}} \to \frac{\Lief}{[\Lies, \Lief]_{\Lie}} \stackrel{\overline{\tau}} \to \Lieq \to 0.
\end{equation}
The above remark shows that there exists at least one \Lie-stem cover associated to a Leibniz algebra having a free presentation $0 \to \Lies \to \Lief \stackrel{\tau}\to \Lieq \to 0$ satisfying $\Lies \subseteq [\Lief, \Lief]_{\Lie}$.

\begin{Co}
All \Lie-stem covers of a given Leibniz algebra $\Lieq$  having a free presentation $0 \to \Lies \to \Lief \stackrel{\tau}\to \Lieq \to 0$ satisfying $\Lies \subseteq [\Lief, \Lief]_{\Lie}$are mutually \Lie-isoclinic.
\end{Co}
{\it Proof.} For two \Lie-stem covers $(g_1)$ and $( g_2)$ of $\Lieq$, $\theta(g_1)$ and $\theta(g_2)$ are isomorphisms by \cite[Proposition 6]{CK}, then $\Ker(\theta(g_1)) = \Ker(\theta(g_2))$ and Proposition \ref{main1} ends the proof. \rdg
\bigskip

\section*{Acknowledgements}

Second author was supported by Ministerio de Economía y Competitividad (Spain) (European FEDER support included), grant MTM2013-43687-P.


\begin{center}

\end{center}

\end{document}